\documentclass[authoryear]{elsarticle}

\usepackage{amssymb,amsmath}
\usepackage{graphics}
\usepackage[arrow, matrix, curve]{xy}
\usepackage{txfonts}
\usepackage{bbm}
\usepackage{verbatim}

\newtheorem{definition}{Definition}[section]

\newtheorem{theorem}{Theorem}
\newtheorem{lemma}[definition]{Lemma}
\newtheorem{corollary}[definition]{Corollary}
\newtheorem{proposition}[definition]{Proposition}

\newtheorem{remark}[definition]{Remark}
\newtheorem{acknowledgment}[definition]{Acknowledgment}

\newproof{pf}{Proof}

\newcommand{\nhat }[1]{\{1,2,\ldots,#1\}}
\newcommand{\ohat }[1]{\{0,1,\ldots,#1\}}

\renewcommand{\Cup}{\bigcup}
\renewcommand{\Cap}{\bigcap}

\newcommand{\N}{\mathbb{N}}

\newcommand{\R}{\mathbb{R}}

\newcommand{\T}{\mathbb{T}}
\newcommand{\Z}{\mathbb{Z}}

\renewcommand{\H}{\mathcal{H}}

\renewcommand{\P}{\mathcal{P}}

\newcommand{\A}{\mathfrak{A}}
\newcommand{\B}{\mathfrak{B}}

\newcommand{\eps}{\varepsilon}
\renewcommand{\phi}{\varphi}
\renewcommand{\rho}{\varrho}
\newcommand{\me}{^{-1}}

\newcommand{\ovl}{\overline}
\newcommand{\xx}{^{(x)}}
\newcommand{\yy}{^{(y)}}
\newcommand{\zz}{^{(z)}}

\makeatletter
\newcommand{\rmnum}[1]{\romannumeral #1}
\newcommand{\Rmnum}[1]{\expandafter\@slowromancap\romannumeral #1@}
\makeatother

\begin{document}
\begin{frontmatter}
\title{Sumset Phenomenon in Countable Amenable Groups}
\author[UNI]{Mathias Beiglb\"ock\fnref{S9612}}
\fntext[S9612]{Supported by the Austrian Science Fund (FWF) under grant S9612.}
\author[OSU]{Vitaly Bergelson\fnref{NSF}}
\fntext[NSF]{Supported by NSF under grant DMS-0600042}
\author[OSU]{Alexander Fish}
\begin{abstract} Jin proved that whenever $A$ and $B$ are sets of positive upper density in $\Z$, $A+B$ is piecewise syndetic. Jin's theorem was subsequently generalized by Jin and Keisler to a certain family of abelian groups, which in particular contains $\Z^d$. Answering a question of Jin and Keisler, we show that this result can be extended to countable amenable groups.
 Moreover we establish that such sumsets (or --- depending on the notation --- ``productsets'') are piecewise Bohr, a result   which for $G=\Z$ was proved by Bergelson, Furstenberg and Weiss.
In the case of an abelian group $G$, we show  that a set is
piecewise Bohr if and only if it contains a sumset of two sets of
positive upper Banach density.
 \end{abstract}
\address[UNI]{Fakult\"at f\"ur Mathematik, Universit\"at Wien\endgraf
Nordbergstra\ss e 15, 1090 Wien, Austria}
\address[OSU]{Department of Mathematics, Ohio State University\endgraf Columbus, OH 43210}
\end{frontmatter}

\section{Introduction}

\subsection{Jin's theorem}
\label{subsec_Jin_thm}

For a set $A \subseteq \Z$, the upper Banach density, $d^{*}(A)$, is defined as
\begin{align}
\label{basic_up_density_formula}
d^{*}(A) = \limsup_{b-a \to \infty} \frac{|A \cap \{a,a+1,\ldots,b\}|}{b-a+1}.
\end{align}
It is well known and not hard to show that if $d^{*}(A)>0$ then the set of differences $A-A = \{a-a' : \, a,a' \in A\}$ is \emph{syndetic}, i.e. has bounded gaps. To see this, one can, for example, argue as follows. First, notice that $n \in A-A$ if and only if $A \cap (A-n) \neq \emptyset$. Second, observe that for any  sequence $(n_i)_{i \in \N} \subset \Z$, the set $A-A$ has to contain an element of the form $n_i-n_j$ for some $i > j$. (This follows from the fact that for some $i > j$ one has to have $(A-n_i) \cap (A-n_j) \neq \emptyset$). Now, if $A-A$ is not syndetic, its complement, $\Z \setminus A$, is \emph{thick}, that is, it contains arbitrarily long intervals. It is easy to see that any thick set in $\Z$ contains a set of differences $D = \{ n_i -n_j, \, i > j\}$ for some sequence $(n_i)_{i \in \N}$. This implies $(A-A) \cap D = \emptyset$ which gives a contradiction.

One cannot expect, of course, that the above fact about the syndeticity of $A-A$ extends to the ``sumset" $A+B = \{a+b : \, a \in A, b \in B\}$ of two arbitrary sets of positive upper Banach density. %For example, one can easily construct thick sets $A,B,C$ such that $A+B \subseteq C$ and $C$ has unbounded gaps.
For example, one can easily construct a thick set $C$ which has unbounded gaps, and such that for some thick sets $A$ and $B$, $A+B \subseteq C$.
 In this case $d^{*}(A) = d^{*}(B) = 1$ but $A+B$ is not syndetic. The following surprising result of Jin shows that, nevertheless, the sumset of any two sets of positive upper Banach density is always \emph{piecewise syndetic}, that is, is the intersection of a syndetic set with a thick set.

\begin{theorem}[\cite{Jin02}]\label{IntegerJin}
Assume that $A,B\subseteq \Z$ have positive upper Banach density. Then there  exist a thick set $C$ and a syndetic set $S$ such that $S\cap C\subseteq A+B$. \end{theorem}

It is not hard to see that not every set of positive upper Banach density is piecewise syndetic. Moreover, one can show that not every set $A$ for which the \emph{density}, $d(A) = \lim_{N \to \infty} \frac{|A \cap \{-N,\ldots,N\}|}{2N+1}$, exists and is positive,
is piecewise syndetic. The following remarks show that any piecewise syndetic set contains a highly structured infinite set of a special type.

Note first that any piecewise syndetic set $S$ in $\Z$ has the property that the union of finitely many shifts of $S$ is a thick set. Now, it is easy to verify that any thick set contains an \emph{IP set}, that is,  a set of the form
$\{x_{n_1}+\ldots+ x_{n_k}:n_1< \ldots < n_k, \, k \in \N\}$, where $(x_n)_{n \in \N}$ is a sequence in $\Z$, which contains infinitely many non-zero elements. Applying Hindman's finite sums theorem, \cite{Hind74}, which states that, for any finite partition of an IP set, one of the cells contains an IP set, we see that any piecewise syndetic set contains a shift of an IP set. On the other hand, one can show that there are sets having density arbitrarily close to $1$ which do not have this property. (This fact was first observed by E.~Strauss, see   \cite[Theorem 2.20]{BBHS06}.)

\subsection{Amenable groups}
\label{subsection_amenable_groups}

It is natural to ask whether  Jin's theorem is valid in a more
general setting where the notion of density can be naturally
formulated. In \cite[Application 2.5]{JiKe03} it is proved that
$A+B$ is piecewise syndetic if $A$ and $B$ are sets which have
positive upper Banach density in $\Z^d$ and recently Jin extended
this result to $\oplus_{d=1}^\infty \Z$ \cite{Ji08}.  Jin and
Keisler \cite[Question 5.2]{JiKe03} ask whether Theorem
\ref{IntegerJin} can be extended to countable amenable groups. In
this paper we answer this Question affirmatively. Before stating our
results we review in this subsection some basic facts about amenable
groups.

A definition of amenability which is convenient for our purposes uses the notion of F\o lner sequence. A sequence $(F_n)_{n \in \N}$ of finite subsets of a countable group $G$ is a (\emph{left})  \emph{F\o lner sequence}  if
\begin{align}
\lim_{n\to \infty} \frac{|gF_n \bigtriangleup F_n|}{|F_n|}=0
\end{align}
for every $g\in G$. Equivalently, $(F_n)_{n\in \N}$ is a F\o lner sequence if for every finite set $K$ and any $\eps>0$ all but finitely many $F_n$ are $(K,\eps)$-invariant in the sense that $|gF_n\bigtriangleup F_n|/|F_n|<\eps$ for all $g\in K$.

 A countable group $G$ is \emph{amenable} if it admits a (left) F\o lner sequence.\footnote{One can show that every amenable group admits also right- and indeed two-sided analogues of left F\o lner sequences. Throughout this paper we deal only with left F\o lner sequences; therefore we will routinely omit the adjective ``left''. }
 The basic
 example of an amenable group is the group of integers, an example of a F\o lner sequence being  an arbitrary
  sequence of intervals $\{a_n,\ldots,b_n\}, n\in \N$ with $b_n - a_n \to \infty$.
  The class of amenable groups is quite rich, and, in particular, contains all solvable groups and is closed under
  the operations of forming directed unions, subgroups and extensions. The basic, but not the only examples of
  non-amenable groups are groups containing the free group on two generators as a subgroup.

Given a set $A$ in an amenable group $G$, denote the relative density of $A$ with respect to a finite set $F$  by
$d_F(A):= \frac{|A\cap F|}{|F|}$. The \emph{upper density of $A$ with respect to a F\o lner sequence $(F_n)_{n \in \N}$} is defined by
\begin{align}
\label{formula3}
\ovl d_{(F_n)}(A):=\limsup_{n\to \infty} d_{F_n}(A),
\end{align}
and we write $d_{(F_n)}(A)$ and call it density with respect to
$(F_n)_{n \in \N}$ if in formula (\ref{formula3}) $\limsup_{n\to
\infty} d_{F_n}(A):=\lim_{n\to \infty} d_{F_n}(A)$. The \emph{upper
Banach density} in amenable groups is defined by
\begin{align}
d^*(A):=\sup\left\{\ovl d_{(F_n)}(A): (F_n)_{n \in \N} \mbox{ is a F\o lner sequence}\right\}.
\end{align}

\begin{remark}
\textnormal{
For $G=\Z$ the above definition differs from original definition of upper Banach density in Subsection \ref{subsec_Jin_thm} (see formula (\ref{basic_up_density_formula})) where the supremum was taken only over intervals instead of arbitrary F\o lner sets. However the two notions are equivalent. 
 For example, this follows from the following general fact which is a simple corollary of Lemma  \ref{ByShiftedSet} below:}

{\it Given a subset $B$ of an amenable group $G$ and any F\o lner
sequence $(F_n)_{n \in \N}$ there is a sequence $(t_n)_{n \in \N}$
such that}
\begin{align}
d^*(B)=d_{(F_nt_n)}(B).
\end{align}
\end{remark}

Given two  sets $A, B$ in a discrete group  $G$ we let $AB=\{ab:a\in A,b\in B\}$.
 A set $S\subseteq G$ is (\emph{left}) \emph{syndetic}  if there is a finite set $F$ such that $FS=G$.
 A set $T\subseteq G$ is called (\emph{right}) \emph{thick}  if for each finite set $F$ there exists some $t\in G$ such that
  $Ft\subseteq T$\footnote{When dealing with non-commutative structures one has at his disposal a ``left/right" choice of notions. For brevity, we just write ``syndetic'' resp.\ ``thick'' for what should rigorously be called 
 ``left syndetic'' resp.\ ``right thick''. The choice of left/right is implicitly present in the definitions of piecewise
  syndetic and piecewise Bohr below.\label{footnote4}}. A set $C\subseteq G$ is \emph{piecewise syndetic} if there exist a thick set $T$ and a  syndetic set $S$ such that $C\supseteq S\cap T$. It is not hard to see that  $C\subseteq G$ is piecewise syndetic if and only if there exists a finite set $K$ such that for each finite set $F$ there is some $t\in G$ such that
$Ft\subseteq KC.$ Piecewise syndetic sets are partition regular: if $C_1\cup\ldots \cup C_r$ is piecewise syndetic,
 then some $C_i,i\in \{1,\ldots,r\}$ is piecewise syndetic. This is not hard to see combinatorially and follows also from the ultrafilter characterization of piecewise syndeticity (cf.\ \cite[Section 4.4]{HiSt98}).

We are now able to state one of the main results of this paper.

\begin{theorem}
%\footnote{Our proof repeated verbatim  extends to the case of locally compact non-compact unimodular amenable groups}
\label{AmenableJin}
Let $G$ be a countable amenable group and
let $A,B\subseteq G$ be such that $d^*(A),d^*(B)>0$. Then $AB$ is piecewise syndetic.
\end{theorem}

%A cautious reader may object that it is not clear that the notion of upper Banach density in amenable groups coincides with the original definition of upper Banach density given for the integers since in this case the supremum was taken only over sequences of intervals instead of arbitrary F\o lner sequences. However the two notions do agree based on the following useful remark which is  a simple consequence of Lemma \ref{ByShiftedSet}.

%\begin{rem}
% Given a subset $B$ of an amenable group $G$ and any F\o lner sequence $(F_n)_{n \in \N}$ there is a sequence $(t_n)_{n \in \N}$ such that
%\begin{align}
%d^*(B)=d_{(F_nt_n)}(B).
%\end{align}
%\end{rem}

%A sequence $(F_n)_{n \in \N}$ of finite subsets of $G$ is a F\o lner sequence if for any finite set $K\subseteq G$ and any $\delta >$  almost all $F_n$ are $(K,\delta)$-invariant.

%Equivalently, $d^*(A)$ is the largest number such that for each $\gamma < d^*(A)$, each $\delta>0$ and each  finite set $ K\subseteq G$ there is $(K, \delta)$-invariant subset $F$ in which $A$ has relative density larger than $\gamma$. \cite{OrWe87, Weis01}

\subsection{Bohr sets. }
\label{subsection_piecewise_Bohr}

The Bohr compactification $bG$ of a countable discrete group $G$ is defined (up to an isomorphism) as the largest  %up to homomorphism)
compact group with the property  that there exists a (not necessarily 1-1) homomorphism $\iota:G \to bG$ which has dense image. While this object exists for very general reasons, it is not always possible to give a useful down-to-earth description of it. Anyway, we will say that  a set $B\subseteq G$  is a \emph{Bohr set} if there exists a non-empty open set $U\subseteq bG$ such that  $B\supseteq \iota \me [U]$.\footnote{The sets $\iota\me [U]$, where $U\subseteq bG$ is open define  the \emph{Bohr-topology }Êon $G$. Hence $B\subseteq G$ is Bohr if and only if it contains a non-empty open set.} If, in the addition, $U$ contains the identity of $bG$ then $B$ will be called Bohr$_0$ set.  If  $G$ is abelian, we can consider the embedding
\begin{eqnarray}
\iota: G &\to& \T^{\hat G}\\
g& \mapsto& (\gamma(g))_{\gamma \in \hat G},
\end{eqnarray}
where $\hat G $ is the dual group of $G$.
Endowed with the product topology, $\T^{\hat G}$ is a compact group, $\overline {\iota [G]}$ is a compact subgroup and it can be shown that it is a ``model" for the Bohr compactification of $G$. This implies that $B\subseteq G$ is a Bohr set if and only if there exist $ \gamma_1,\ldots , \gamma_n\in \hat G$ and an open set $U\subseteq \T^n$ such that $\{g\in G:\gamma_1(g),\ldots , \gamma_n(g)\in U\}$ is non-empty and contained in $B$.

Call a set $A\subseteq G$ \emph{ piecewise Bohr } if it is the intersection of a Bohr set and a thick set. Since every Bohr set is syndetic, piecewise Bohr sets are piecewise syndetic.

 By \cite[Theorem 4.3]{BeFW06} there exists a syndetic set of integers which is not piecewise Bohr. Note that this also implies that there exists a partition of the integers into finitely many cells none of which is piecewise Bohr.
 % Neil Hindman has a relatively simple combinatorial counterexample which was never published, we should mention it in some appropriate way.\marginpar{Mathias, let us kill it}

Given a Bohr set $B$ there exist a Bohr$_0$ set $B_0$ and a Bohr set $B_1$ such that  $B\supseteq B_0B_1$. This is a trivial consequence of the fact that the Bohr-topology is a group topology on $G$. Also, given a thick set $T$, it is not difficult to see that there exist thick sets $T_0$ and $T_1$ such that $T\supseteq T_0T_1$ provided that $G$ is abelian. (See Lemma \ref{ThickIsDivisible} below.)  It follows that for every piecewise Bohr set $A$ there exist piecewise Bohr sets $A_0,A_1$ such that $A_0A_1\subseteq A$.  In particular every piecewise Bohr set contains the product of two sets of positive upper Banach density. This puts an upper bound on the amount of structure which can be expected in the productset of two sets of positive upper Banach density. Somewhat surprisingly, it is in fact always possible to get this much:

\begin{theorem}\label{AmenableBohrness}
Let $G$ be a countable amenable group and assume that $A,B\subseteq
G$ have positive upper Banach density. Then $AB$ is piecewise Bohr.
\end{theorem}

In the case $G=\Z$, Theorem \ref{AmenableBohrness} is proved in \cite{BeFW06}.

Summarizing the above discussion, we have the following characterization of sumsets in the abelian case.

\begin{theorem}
\label{cor_pwBohr_char}
Let $(G,+)$ be a countable abelian group and let $C\subseteq G$. Then $C$ is piecewise Bohr if and only if there exist sets $A, B$ of positive upper Banach density such that $A+B \subseteq C$.
\end{theorem}

 We will show in Section \ref{AbelianSection} that Theorem \ref{cor_pwBohr_char} does not extend to the non-commutative setup.

%We will see in Proposition \ref{ThickIsNotDivisible} below that in
%the discrete Heisenberg group there exists a piecewise Bohr set $C$
%(in fact $C$ will be thick)  which does not contain the product set
%$AB$ where sets $A$ and $B$ have positive upper Banach density
%((left, left), (right, right), (left, right) or (right, left)).

\subsection{Organization of the paper}

In Section \ref{IntegerSection} we provide a simple proof of Jin's
Theorem for $G=\Z$.  In  Section  \ref{EasyAmenable} we explain how
this proof can be modified to extend Jin's result to the amenable
setting (Theorem \ref{AmenableJin}). 
The results in Section
\ref{DDsection} allow us to give yet another proof of Theorem \ref{AmenableJin} and will also be utilized in Section \ref{BohrSection} in the proof of Theorem
\ref{AmenableBohrness}. 
 Finally, in Section \ref{AbelianSection} we
prove Theorem \ref{cor_pwBohr_char} and provide an example which
demonstrates that  Theorem \ref{cor_pwBohr_char} does not extend to
the non-commutative setup.

Throughout this paper, $G$ will denote a countable discrete amenable group.
We call $(X,\B,\mu)$ a Borel probability space if $(X,\B)$ is a measurable space isomorphic to the unit interval equipped with the $\sigma$-algebra of Borel sets and $\mu$ is  a Borel probability measure on $(X,\B)$. If $(X,\B, \mu)$ is a Borel probability space and $T:X\to X$ is an invertible measure preserving transformation, $(X,\B,\mu,T)$ will be called a measure preserving system.

\section{Jin's theorem in the integers}\label{IntegerSection}
Jin's original proof of Theorem \ref{IntegerJin} in \cite{Jin02} utilized non-standard analysis. Jin also  provided a purely combinatorial proof of Theorem \ref{IntegerJin} (\cite{Jin04}). The purpose of this ``warm-up" section is to give another proof of Theorem \ref{IntegerJin}.
While our proof is shorter than the original one, most of the ideas we use can be found, at least implicitly, in Jin's work.

Our proof of Jin's theorem will be based on the following two lemmas:

\begin{lemma}\label{ThickLemma} Assume that $A,B$ are sets of integers such that $d^*(A)+d^*(B)>1$. Then $d^*(A+B)=1$, i.e.\ $A+B$ is thick.
\end{lemma}

\begin{lemma}\label{GrowthLemma}\hspace{-2mm}\footnote{
Lemma  \ref{GrowthLemma}  is originally due to Neil Hindman, see  \cite[Theorem 3.8]{Hind82}. The combinatorial proof given subsequently is based on the same idea as Hindman's proof.}\hspace{0mm} If $A$ is a set of integers then $\sup_{k\geq 0} d^*(\{-k,\ldots, k\}+A)$ is either $0$ or $1$.
\end{lemma}

Taking Lemmas \ref{ThickLemma} and \ref{GrowthLemma} for granted, Theorem \ref{IntegerJin} is almost trivial: By Lemma \ref{GrowthLemma} there is some integer $k$ such that $d^*(\{-k,\ldots, k\}+A )+d^*(B)>1$. Hence by Lemma \ref{ThickLemma}, $\{-k,\ldots, k\}+A+B$ is thick. Thus $A+B$ is piecewise syndetic.

\medskip
 Recall that for a finite interval $I\subseteq \Z$ and a set $A\subseteq \Z$,
$d_I(A)=\frac {|I\cap A|}{|I|}$
denotes the relative density of $A$ with respect to $I$.
\newproof{potTL}{Proof of Lemma \ref{ThickLemma}}
\begin{potTL}
%Since $d^*(-A)=d^*(A)$, it is sufficient to prove that $-A+B$ is thick.
Note that if $J\subseteq \Z$ is \emph{any} non-empty interval and $d^*(B)>\beta$, then there exists  $t\in \Z$ such that $d_{J+t}(B)>\beta$.

Pick $\alpha, \beta>0$ such that $d^*(A)>\alpha, d^*(B)>\beta, \alpha+\beta =1$ and fix $n\in \N$. We have to prove that $A+B$ contains a shifted copy of $\ohat n$. Loosely speaking, long enough intervals are almost invariant with respect to shifts by elements of $\ohat n$. In particular there exists an interval $I$ such that $
d_I(-x+A)>\alpha$ for all $x\in \ohat n$.

Apply the above observation to the interval $J=-I$ and
pick  some integer $t\in \Z$ such that $d_{(-I)+t}(B)=d_{-I}(B-t)>\beta$.
Let $x\in \ohat n$.  Since $\alpha + \beta =1$,
\begin{align}
d_{-I}(-A+x)+ d_{-I}(B-t)>1  \ & \Rightarrow \ (-A+x)\cap (B-t)\neq \emptyset \\
&  \Rightarrow \ x+t\in A+B.
\end{align}
Since $x$ was arbitrary, we have $\ohat n +t \subseteq A+B$ as required.
\qed\end{potTL}

\medskip

We will give two proofs of Lemma \ref{GrowthLemma}. The first one is based on an elementary  combinatorial argument, the second one involves more abstract concepts and gives a rigorous meaning to the
intuitive fact expressed by Lemma  \ref{GrowthLemma} that the system $(\Z, \P(\Z), n\mapsto n+1, d^*)$ is ``ergodic".

\newproof{potCGL}{Combinatorial proof of Lemma \ref{GrowthLemma}}
\begin{potCGL}
We will show that for any set $A\subseteq \mathbb Z$ with $d^*(A)>0$ one has $\sup_{n\geq 0} d^*(A+ \{-n,\ldots, n\})= 1$. Assume by way of contradiction that $d^*(A)>0$, but $\sup_{n\geq 0} d^*(A+ \{-n,\ldots,n\})= \gamma < 1$. Pick $\eps>0$ such that $(\gamma+\eps)^2 <\gamma-\eps$.
For $n$ large enough, $d^*(A+ \{-n,\ldots, n\})>\gamma-\eps$.
Hence, replacing $A$ by $A+\{-n,\ldots, n\}$ if necessary, we may assume that $d^*(A)>\gamma-\eps$.

Fix $k\in \N$ such that $d_I(A)<\gamma+\eps $  for any interval $I\subseteq \Z$ of length $k$. Then pick an interval $J$ such that the following conditions are satisfied:
\begin{enumerate}
\item The length of $J$ is $m\cdot k$ for some positive integer $m$.
\item $d_J(A+\{-k, \ldots, k\})< \gamma+\eps$.
\item $d_J(A)>\gamma-\eps$.
\end{enumerate}
Partition $J$ into intervals $I_1,I_2,\ldots, I_m$ of length $k$. Assume  that $A$ intersects more than  $m\cdot (\gamma+\eps)$ of these intervals. Then $A+\{-k, \ldots, k\}$ \emph{covers} more than $m\cdot (\gamma+\eps)$ of these intervals, hence $d_J(A+\{-k, \ldots, k\})$ exceeds $m\cdot (\gamma+\eps)/m=\gamma+\eps$, contradiction. Thus $A$ intersects at most $m\cdot (\gamma+\eps)$ of the intervals $I_j, j\in \nhat m$. Since the relative density of $A$ in a length $k$ interval is bounded by $\gamma+\eps$ this yields
\begin{align}
d_J(A)\leq (\gamma+\eps)\cdot m\cdot (\gamma+\eps)/m=(\gamma+\eps)^2
\end{align}
which contradicts $(\gamma+\eps)^2 <\gamma-\eps$.
 \qed\end{potCGL}
 
 \medskip

Our second proof of Lemma \ref{GrowthLemma} is based on the following version of Furstenberg's correspondence principle.
\begin{proposition}\label{CorrespondencePrinciple}
Assume that $A\subseteq \Z$ has positive upper density. Then there exist an \emph{ergodic} measure preserving system $(X,\B,\mu,T)$ and a measurable set $B\subseteq X$ such that
\begin{align}
d^*(A)& = \mu(B) \mbox { and}\\
\label{CPBound} d^*(A-n_1\cup \ldots \cup A-n_k ) & \geq \mu(T^{-n_1}B\cup \ldots \cup T^{-n_k} B)
\end{align}
for all $n_1,\ldots, n_k\in \Z$.
\end{proposition}

Proposition \ref{CorrespondencePrinciple} differs from the more familiar forms of Furstenberg's correspondence principle (see \cite[Theorem 1.1]{Ber85}) in that we use unions instead of intersections and in that we require that
$(X,\B,\mu,T)$ to be ergodic.
One can easily verify that due to the algebraic nature of Furstenberg's correspondence principle, virtually any known proof (see, in particular, the proofs in \cite{Ber85, BeMc98})  is equally valid for unions.
That the system can be
chosen to be ergodic follows from  \cite[Proposition 3.9]{Furs81}.

%Proposition \ref{CorrespondencePrinciple} is almost identical to  \cite[Proposition 3.1]{BeHK05} with the only difference that in  (\ref{CPBound})  a bound holds for \emph{unions} instead of \emph{intersections}. However, the proof of Proposition \ref{CorrespondencePrinciple} is precisely the same, hence we will not state it here.
\newproof{potDGL}{``Dynamical'' proof of Lemma \ref{GrowthLemma}}
\begin{potDGL}
Assume that $d^*(A)>0$ and choose $(X,\B,\mu,T)$ and $B\subseteq X$ according to Proposition \ref{CorrespondencePrinciple}. Since $T$ is ergodic, $$  \sup _{k\geq0} d^*\big(\{-k,\ldots, k\}+A\big)\geq \sup _{k\geq0} %\mu\Big(\bigcup_{k=-n}^n T^{-k} B\Big) \geq 
\mu\big(T^{-k}B\cup\ldots \cup T^k B\big) = 
\mu\Big(\bigcup_{k\in \Z} T^{-k}B\Big)=1. \mbox{\qed}$$
\end{potDGL}

\begin{remark}
For the usual (upper) density the statement of Lemma \ref{GrowthLemma} is not true. For example, let
 $B=\bigcup_{n\in \N} \{n^2, n^2+1,\ldots, n^2+n\}$. Then for $ A = B \cup (-B)$ we have
\begin{align} d(A)=\lim_{N \to \infty} \frac{|A \cap [-N,\ldots,N]|}{2N+1} =1/2=\sup_{k\geq 0}\overline d(\{-k,\ldots,k\}+A).\end{align}
\end{remark}

However,  it follows from the proof of Lemma \ref{AmenableGrowthLemma}, that if $(F_n)_{n \in \N}$ is a  F\o lner sequence  which satisfies $d_{(F_n)}(A)=d^*(A)>0$), then we have
\begin{align}
 \sup_{k\geq 0}\overline d_{(F_n)}(\{-k,\ldots,k\}+A)=1.
 \end{align}

%\medskip

\section{Jin's theorem in countable amenable groups}\label{EasyAmenable}

In this section we demonstrate that (with some work) the proof of Jin's theorem which was given in the previous section generalizes to the amenable setting.
%Subsequently we will explain in which way the arguments can be adapted.
The proof of the  general ``amenable" statement is  based on the following auxiliary results. (cf.\ Lemmas \ref{ThickLemma}, \ref{GrowthLemma})% which are natural amenable analogues of the corresponding statements in the integer setting.

\begin{lemma}\label{AmenableThickLemma}
Let $G$ be an amenable group and assume that $A,B\subseteq G$, $d^*(A)+d^*(B)>1$. Then $A B$ is thick.
\end{lemma}

\begin{lemma}\label{AmenableGrowthLemma}
Let $G$ be a countable amenable group and let $A\subseteq G$. Then $\sup \{d^*(KA): K\subseteq G, K \mbox{ is finite}\}$ is either $0$ or $1$.
\end{lemma}

Note first, that in complete analogy with the integer setting, Lemma \ref{AmenableThickLemma} and Lemma \ref{AmenableGrowthLemma} imply that if $d^*(A),d^*(B)>0$, then there exists a finite set $K$ such that $KAB$ is thick,
which, in turn, implies that  $AB$ is  piecewise syndetic.

The following simple fact is needed in the proof of Lemma \ref{AmenableThickLemma} (and will also be utilized in the next section for the proof of Lemma \ref{Hurray}.).
\begin{lemma}\label{ByShiftedSet}
Let $B, K\subseteq G$, $K$ finite and $\beta < d^*(B)$. Then there exists some $t\in G$ such that
\begin{align}\label{ShiftedSet}
d_{Kt}(B)=\frac{|B\cap Kt|}{|K|}\geq \beta.
\end{align}
\end{lemma}

\begin{pf}{Proof.}
Pick a  F\o lner set $F$ such that $|B\cap gF| /|F|\geq \beta$ for each $g\in K$. Then
\begin{align}
\sum_{t\in F} |B\cap Kt|=|\{(g,t)\in K\times F:gt\in B\}|= \sum_{g\in K}|B\cap gF|\geq |K| \cdot  |F| \cdot \beta.
\end{align} Dividing by $|K| \cdot  |F|$ we see that (\ref{ShiftedSet}) holds for some $t\in F$.
\qed\end{pf}

\medskip

\newproof{potATL}{Proof of Lemma \ref{AmenableThickLemma}}
\begin{potATL}
To obtain  Lemma \ref{AmenableThickLemma}, one just has to rewrite the proof of Lemma \ref{ThickLemma} in terms of F\o lner sequences. The only part which needs  justification is that if $d^*(B)>\beta $ and $F\subseteq G$ is a finite set, then there is some $t\in G$ such that $d_{F\me t}(B)>\beta$. This was proved in Lemma \ref{ByShiftedSet}.
\qed\end{potATL}

\medskip

 Lemma \ref{AmenableGrowthLemma} can be proved in a variety of ways. First, it is possible to prove  an appropriate  version of Furstenberg's correspondence principle for amenable groups (for instance, one can combine the proof of
correspondence principle given in \cite[Theorem 2.1]{BeMc98} or in \cite[Theorem 6.4.17]{Be00} with the amenable analogue of \cite[Proposition
3.9]{Furs81}) which then immediately  gives the desired  result as in the dynamical proof of Lemma \ref{GrowthLemma}.
 %It seems plausible that many experts in ergodic are well aware how the appropriate correspondence principle is obtained and that a rigorous proof would not appeal to non specialist in ergodic theory, hence we refrain from giving details here.

Second, one also can prove  Lemma \ref{AmenableGrowthLemma} via an appropriate generalization of the combinatorial proof of Lemma \ref{GrowthLemma}. There we employed  the fact that intervals \emph{tile} the integers. In general, a set $T$ in a countable group $G$ is a \emph{tile}  if there exists a set $S\subseteq G$ such that $\{Ts:s\in S\}$ is a partition of $G$. %Given a finite subset $K$ of an abelian group $G$ and $\eps>0$, it is easy to construct a $(K, \eps)$-invariant tile of $G$. The subgroup $G'=\langle K\rangle$ is of the form $\Z^k\times G_0$, where $k\geq 0$ and $G_0$ is a finite abelian group. If $n$ is chosen sufficiently large, $T:=\{0,\ldots n\}^k\times G_0$, $T$ is $(K,\eps)$ invariant. Trivially $T$ is a tile of $G'$ and since $G'$ is a subgroup of $G$, $T$ is also a tile of $G$. Thus every abelian groups has a F\o lner sequence all of whose elements are tiles of $G$. This observation is sufficient to extend Lemma \ref{GrowthLemma} to the general abelian setting.
The group
 $G$ is called \emph{monotilable} if it admits a F\o lner sequence consisting of tiles and in this case the proof of Lemma \ref{GrowthLemma} can be adapted fairly naturally. Having the construction of F\o lner sequences in the abelian setting in mind, it is easy to see that every countable abelian group is monotilable and it is shown in \cite{Weis01} that much more general classes of amenable groups share this property.  While it is not known whether all amenable groups are monotilable, they do admit so called quasi-tilings (see \cite{OrWe87}). Those still do allow to push the proof of Lemma \ref{GrowthLemma} to the desired generality, but the details become unpleasantly technical.

Since Lemma \ref{AmenableGrowthLemma} is crucial for  a generalization of  Jin's theorem to the amenable case, we will give here a self contained proof. While the argument is more involved than that used in the combinatorial proof of Lemma \ref{GrowthLemma}, it is still entirely elementary.

\newproof{potAGL}{Proof of Lemma \ref{AmenableGrowthLemma}}
\begin{potAGL} It is sufficient to consider the case $d^*(A)>0$.
Pick   a F\o lner sequence $(F_n)_{n \in \N}$ such that $d_{(F_n)}(A)=\alpha>0$ and $d_{(F_n)}(KA)$ exists for each finite $K\subseteq G$. Let $\beta = \sup\{d_{(F_n)}(KA): K\subseteq G, \mbox{$K$ finite}\}$.
We claim that after passing, if necessary, to a subsequence of $(F_n)$, there exists a F\o lner sequence $(G_n)_{n \in \N}, G_n\subseteq F_n$ such that the following hold true:
\begin{enumerate}
\item $\lim_{n\to \infty } |G_n|/ |F_n|=\beta$.
\item $d_{(G_n)}(HA)  = d_{(F_n)}(HA)\frac 1 \beta $ for any finite set $H\subseteq G$.
\end{enumerate}
A particular consequence of (2) is that  $\sup\{d_{(G_n)}(KA): K\subseteq G, \mbox{$K$ finite}\}=\beta/\beta= 1.$

Fix a sequence $(K_n)_{n \in \N}$ of finite subsets of $G$  such that  $K_n K_n \subseteq K_{n+1}$, $K_n \uparrow G$ and each $K_n$ contains the identity of $G$. Passing to a subsequences once more, we can assume that
\begin{align}
d_{F_m}(K_nA)\in (\beta-1/n, \beta+1/m) \mbox{ for all $m\geq  n$,}
\end{align} and that each $F_n$ is $(K_n, 1/n)$-invariant. Set $G_n:=K_{n-1} (A\cap F_n)$. Note that
\begin{align}
G_n&\subseteq K_{n-1}A \cap K_{n-1} F_n \approx K_{n-1} A \cap F_n\\
G_n&\supseteq K_{n-1}A \cap \bigcap_{k\in K_{n-1}} k F_n \approx K_{n-1} A \cap F_n,
\end{align}
since $F_n$ is assumed to be almost invariant with respect to $K_{n-1}$. In particular $|G_n|/|F_n|\approx \beta$. Next we show that $(G_n)_{n \in \N}$ is a F\o lner sequence. To this end, fix $n\in \N$ and $t\in K_{n-1}$.  We have
\begin{align}
& \frac{|tG_n\setminus G_n |}{|G_n|}  \approx\frac{|(tK_{n-1}A \cap F_n)\setminus G_n |}{|G_n|}\subseteq \\
 \subseteq & \frac{|(K_{n}A \cap F_n)\setminus G_n |}{|G_n|} \approx \frac{|(K_{n}A \cap F_n)\setminus (K_{n-1}A \cap F_n) |}{|F_n|} \frac1\beta = \\
 = & \left(\frac{|K_{n}A \cap F_n|}{|F_n|}-\frac{|K_{n-1}A \cap F_n|}{|F_n|}   \right)\frac1\beta = \big(d_{F_n}(K_n A)-d_{F_n}(K_{n-1}A)\big) \frac 1\beta \to 0.
\end{align}
Finally we have
\begin{align}
d_{(G_n)}(HA)=&\lim_{n\to\infty} \frac{| HA \cap (K_{n-1}A \cap F_n)|}{|K_{n-1}A \cap F_n|}\\
=& \lim_{n\to\infty} \frac{| (HA \cap K_{n-1}A) \cap F_n|}{\beta | F_n|} = \frac1\beta d_{(F_n)}(HA),
\end{align}
which gives us (\rmnum{2}).
\qed\end{potAGL}

\medskip

\section{Finer structure of productsets.}\label{DDsection}

%The following proposition (which is the main result of this section) shows  that the product of two sets of positive upper Banach density contains a certain specific structure. 
The following proposition (which is the main result of this section) shows  that the product of two sets of positive upper Banach density contains translations arbitrarily large pieces of the product of a ``large set'' with its inverse. 
(This fact  will be
utilized in the proof of Theorem \ref{AmenableBohrness} in the next
section.)

\begin{proposition}\label{PiecewiseDD} Let $G$ be a countable amenable group and
let $A,B\subseteq G$ be such that $d^*(A),d^*(B)>0$. Then there exists a set $D\subseteq G$ with $d^*(D)>0$ such that for each finite set $H\subseteq G$, there is some $t_H$ such that
\begin{align}
(H\cap DD\me) t_H\subseteq AB.
\end{align}
\end{proposition}
Using Lindenstrauss' pointwise ergodic theorem \cite{Lind01} it is possible to show that for any set $D$ which has positive upper Banach density and for any F\o lner sequence $(F_n)_{n \in \N}$ to which the pointwise ergodic theorem applies, there exists a set $E$ such that $d_{(F_n)}(E)= d^*(D)$ and $EE\me \subseteq DD\me$. Hence it is possible to give a somewhat stronger formulation of Proposition \ref{PiecewiseDD}.

%Note the following immediate consequence of Lemma \ref{ByShiftedSet} which provides some basic information about the behavior    of upper Banach density in amenable groups.

%\begin{corollary}\label{BanachByShiftedFolner}
%Given $B\subseteq G$ and a F\o lner sequence  $(F_n)_{n \in \N}$ there exists a sequence $(t_n)_{n \in \N}$ in $G$ such that \begin{align}d_{(F_nt_n)}(B)=d^*(B).\end{align} \end{corollary}

Before proving Proposition \ref{PiecewiseDD} we formulate and prove
a few auxiliary results.

\begin{lemma}\label{Hurray}
Let $A_0, B\subseteq G$, $A_0$ finite and $\beta < d^*(B)$. There exist $C\subseteq A_0$ and $ t\in G$ such that $CC\me t\subseteq  A_0 B$ and $ |C|\geq \beta |A_0|.$
\end{lemma}

\begin{pf}
Applying Lemma \ref{ByShiftedSet} to $A_0\me$ we find $t$ such that \begin{align} \beta |A_0| \leq |A_0\me t \cap B|= |A_0\cap(Bt\me)\me|.\end{align} And for all $x, y \in C:= A_0\cap(Bt\me)\me$ we have
$xy\me\in A_0(Bt\me)$.
\qed\end{pf}

\medskip

While the formulation of Lemma \ref{Hurray} appears to be somewhat technical, it allows to
show that $AB$ contains arbitrary large sets of the form $C_nC_n\me t_n $.
%squeeze increasingly large sets of the form $C_nC_n\me t_n $ inside $AB$.
%Hence in a way it takes us a great deal closer   towards Proposition \ref{PiecewiseDD}.
The remaining ingredient in the proof of Proposition \ref{PiecewiseDD} is the following statement.
%which close the gap we need the following lemma.

\begin{lemma}\label{NiceLemma}
Let $(F_n)_{n \in \N}$, $(G_n)_{n \in \N}$ be F\o lner sequences, let $C_n\subseteq F_n$ and set  $\gamma:=\limsup d_{F_n}(C_n)$. Then there exists a set $D$ such that the following hold.
\begin{enumerate}
\item $\ovl d_{(G_n)}(D) =\gamma$.
\item For each finite set $D_0\subseteq D$ there exist $c\in G$ and $n\in \N$ such that $D_0c\subseteq C_n$.
\end{enumerate}
%In particular, given a sequence $(t_n)_{n \in \N}$, we can choose for each finite set $K\subseteq D$ some $t_K\in G$ such that \begin{align}\label{KKtK} \Cup_{K\subseteq D, |K|<\infty} KK\me t_K\subseteq \Cup_n C_nC_n\me t_n. \end{align}
\end{lemma}

The proof of Lemma \ref{NiceLemma} relies on the following Fubini-type Lemma. % which is a somewhat extended version of \cite[Theorem 1.1]{Berg85}.
 \begin{lemma}\label{VitalyFubini}
Let $(X, \A, m)$ be some space equipped with a finitely additive measure, assume that $(A_g)_{g\in G}$ is a sequence of sets in $\A$ such that $m(A_g)\geq \gamma$ for all $g\in G$ and let $(G_n)_{n \in \N}$ be a F\o lner sequence. Then there exists a set $D$ such that $\ovl d_{(G_n)}(D)\geq \gamma$ and $m(\Cap _ {t\in D_0} A_t)>0$ for every finite set $D_0\subseteq D$.
\end{lemma}
% \begin{pf*}{Sketch of Proof.}
Lemma \ref{VitalyFubini} is essentially \cite[Lemma 5.10]{Berg06}, the only difference being that here we only require that $m$ is finitely additive. The following argument shows that the case of finitely additive measures follows from the $\sigma$-additive setup. Indeed,  set $Y:=\{0,1\}^\N$, let $B_n=\{(x_k)_{k \in \N}\in Y: x_n=1\}$ for $n\in\N$ and
put
\begin{align}
\mu_0\Big(\Cap_{k\in S} B_k\cap \Cap_{n\in T} (Y\setminus B_n)\Big):=m\Big(\Cap_{k\in S} A_k\cap \Cap_{n\in T}(Y\setminus A_n)\Big)
\end{align}
for finite sets $S,T\subseteq \N$. Then $\mu_0$ naturally extends to a $\sigma$-additive Borel probability measure $\mu$ on  $Y$ and it is sufficient to prove Lemma \ref{VitalyFubini} for the sets $B_1,B_2,\ldots $ in $(Y, \B,\mu)$.

% \cite[Theorem 1.1]{Berg85} asserts that if $B_n,n\in\N$ are sets in a probability measure space $(Y,\B,\mu)$ all having measure at least $\gamma$, then there exists a set $D$ which has upper density at least $\gamma$ such that $\Cap_{t\in D_0} B_t$ has positive measure for each finite set $D_0\subseteq D$.
%Lemma \ref{VitalyFubini} is almost the same assertion as \cite[Theorem 1.1]{Berg85}, the sole differences being that the original result deals with usual upper density in $\N$ and refers to a countably additive measure space.
%An inspection of the argument given in \cite{Berg85} reveals that passing from intervals of integers to a sequence $(G_n)_{n \in \N}$ of finite sets does not affect the proof.  Neither does it significantly alter the problem to drop the assumption of $\sigma$-additivity  since it is possible to consider a countably additive model of $(X,\A,m)$. Indeed set $Y:=\{0,1\}^\N$, let $B_n=\{(x_k)_{k \in \N}\in Y: x_n=1\}$ for $n\in\N$ and put
%\begin{align}
%\mu_0\Big(\Cap_{k\in S} B_k\cap \Cap_{n\in T} (Y\setminus B_n)\Big):=m\Big(\Cap_{k\in S} A_k\cap \Cap_{n\in T}(Y\setminus A_n)\Big)
%\end{align}
%for finite sets $S,T\subseteq \N$. Then $\mu_0$ naturally extends to a Borel probability measure $\mu$ on  $Y$ and it is sufficient to prove Lemma \ref{VitalyFubini} for the sets $B_1,B_2,\ldots $ in $(Y, \B,\mu)$.

%\qed \end{pf*}

\newproof{potNL}{Proof of Lemma \ref{NiceLemma}}
\begin{potNL}
Passing to a subsequence if necessary, we can assume that $\gamma=\lim d_{F_n}(C_n)$ exists. Consider $C:=\Cup_n C_n\times \{n\}\subseteq G\times \N=: X$. Let $\A$ be the algebra of subsets of $X$ generated by all sets of the form $gC:=\Cup_n (gC_n)\times \{n\}, g\in G$. Since $\A$ is countable, we can pick a sequence $k_1<k_2< \ldots$ in $\N$ such that
\begin{align}
m(A)=\lim_{k\to \infty} \frac{|A\cap (F_{n_k}\times \{n_k\})|}{|F_{n_k}|}
\end{align}
exists for all $A\in \A$. Note that
\begin{align}
m(gC)=\lim_{k\to \infty} \frac{|gC_{n_k} \cap F_{n_k}|}{|F_{n_k}|}=\gamma
\end{align}
for all $g\in G$.
 Let $D$ be the ``outcome" of applying Lemma \ref{VitalyFubini} to the space $(X, \A,m)$ and the sets $g\me C, g\in G$. Given a finite set ${D_0}\subseteq D$, $m\left( \Cap_{g\in {D_0}}g\me C\right)>0$. Hence for $k$ large enough, $\Cap_{g\in {D_0}}g\me C_{n_k}$ has positive relative density with respect to $F_{n_k}$, so pick $c\in \Cap_{g\in {D_0}}g\me C_{n_k}$. Then ${D_0}c \subseteq C_{n_k}$ as required.
\qed\end{potNL}

\medskip

We are now in the position to prove the main result of this section.

\newproof{potDD}{Proof of Proposition \ref{PiecewiseDD}} 
\begin{potDD} Pick a F\o lner sequence $(F_n)_{n \in \N}$ and $\alpha >0$  such that $d_{F_n}(A)\geq \alpha>0$ for all $n\in\N$. Pick $\beta>0$ such that $d^*(B)>\beta$. Applying Lemma \ref{Hurray} to the sets $A_n:=A\cap F_n$, we find sequences $(C_n)_{n \in \N}$ and $(t_n)_{n \in \N}$ such that $\bigcup_{k=1}^\infty C_kC_k\me t_k\subseteq AB$ and $d_{F_n}(C_n)\geq \alpha \beta>0$ for each $n\in \N$.
Pick a set $D$ guaranteed by Lemma \ref{NiceLemma}. Given an arbitrary finite set $H\subseteq G$, there is a finite set $D_0\subseteq D  $ such that $H\cap DD\me \subseteq D_0D_0\me$. By Lemma \ref{NiceLemma}, there exist $c\in G$ and $n\in \N$ such that $D_0c\subseteq C_n$. Hence 
\begin{align}
(DD\me \cap H)t_n\subseteq D_0D_0\me t_n=D_0c (D_0c)\me t_n\subseteq \bigcup_{k=1}^\infty C_kC_k\me t_k\subseteq AB.
\end{align}
\qed \end{potDD}

\medskip

In the next section we will use Proposition \ref{PiecewiseDD} together with Lemma \ref{DifferencesContainReturnTimes} to prove that $AB$ is piecewise Bohr if $d^*(A),d^*(B)>0$.

\begin{lemma}\label{DifferencesContainReturnTimes}
Let $A\subseteq G$ and assume that $d^*(A)>0$. Then there exist a Borel probability space $(X,\B,\mu)$, a measure preserving $G$ action $(T_g)_{g\in G}$ on $X$ and a set $B\subseteq X, \mu(B)=d^*(A)$ such that
\begin{align}
\{g\in G: \mu(T_{g}\me B\cap B)>0 \} \subseteq AA\me.
\end{align}
In particular $AA\me $ is syndetic.
\end{lemma}

In a certain sense Lemma \ref{DifferencesContainReturnTimes} can be reversed. Indeed, using the ergodic theorem it is not difficult to see that for  any set $R$ of return times there exists a set $A$ of positive upper Banach  density such that $AA\me \subseteq R$.

\medskip

We will derive  Lemma \ref{DifferencesContainReturnTimes} from the following amenable version of Furstenberg's correspondence principle (see for instance \cite[Theorem 2.1]{BeMc98}, \cite[Theorem 6.4.17]{Be00}).

\begin{lemma}
\label{AmenbleCorrespondencePrinciple}
Let $G$ be an amenable group and assume that $A\subseteq G$. Then there exist a Borel probability space $(X,\B,\mu)$,  a measure preserving $G$ action $(T_g)_{g\in G}$ on $X$ and set $B\subseteq X, \mu(B)=d^*(A)$ such that
\begin{align}
 \mu(T_{g_1}\me B\cap\ldots\cap  T_{g_n}\me B)\leq d^*(g_1\me A\cap \ldots\cap g_n\me A)
\end{align}
for all $g_1,\ldots, g_n\in G$.
\end{lemma}

\newproof{potDC}{Proof of Lemma \ref{DifferencesContainReturnTimes}}
\begin{potDC}
Let $(X, \B, \mu), (T_g)_{g\in G}$ and $B$ be as in  Proposition \ref{AmenbleCorrespondencePrinciple}. Then
\begin{align}
AA\me \supseteq \{g:d^*(g\me A\cap A)>0\}\supseteq \{g:\mu (T_g\me B\cap B)>0\}=:S.
\end{align}
Set $Y:= \bigcup_{g\in G} T_g\me B$.  Pick a finite set $K\subseteq G$ such that $\mu(\bigcup _{g\in K} T_g\me B)+\mu(B)>\mu(Y)$. Fix $h\in G$. Then $\mu(\bigcup _{g\in K} T_g\me B\cap T_h\me B)>0$.
Hence for some $g\in K$ we have $\mu(T_{gh}\me B\cap B)>0$. Equivalently $gh\in S=\{f\in G: \mu(T_f\me B\cap B)>0\}$. Since $h\in G$ was arbitrary, $G=K\me S$, so $AA\me $ is indeed syndetic.
\qed\end{potDC}

\medskip

We conclude this section with showing how Proposition \ref{PiecewiseDD} offers yet another way to establish Theorem \ref{AmenableJin}. If $A,B$ have positive upper Banach density, we may choose a set $D$ of positive upper Banach density such that $AB$ contains shifts of arbitrary finite portions of $S=DD\me$. By Lemma \ref{DifferencesContainReturnTimes} the set $S$ is syndetic and hence also piecewise syndetic. 
Thus piecewise syndeticity of $AB$ follows from the following natural property of piecewise syndetic sets. 

\begin{lemma} \label{PSRegularity}
Let $G$ be a group, $S,T\subseteq G$ and assume that $S\subseteq G$ is piecewise syndetic and that for each finite set $H\subseteq G$ there is some $t_H\in G$ such that
\begin{align}
(H\cap S)t_H \subseteq T.
\end{align}
Then $T$ is piecewise syndetic as well.
\end{lemma}

\begin{pf}
 Pick  a finite set $K\subseteq G$ such that $KS $ is thick. Given an arbitrary finite set $F\subseteq G$, there is some $f\in G$ such that $Ff\subseteq KS$. Choose a finite set $H$ such that $F\subseteq K(S \cap H)$. Since $(S \cap H)\subseteq T t_H\me$, we have $Ff\subseteq KTt_H\me.$ As $H$ was arbitrary, $KT$ is thick.
\qed\end{pf}

\medskip

%To prove Corollary \ref{AmenableJin} from Theorem \ref{PiecewiseDD} it is sufficient to  have that $DD\me$ is syndetic (or just piecewise syndetic). The proof that we have given for this statement does not use $\sigma$-additivity, hence it would be possible to avoid using the Correspondence Principle. However we will need it in the next section.

\section{Bohr sets and almost periodic functions}\label{BohrSection}

Consider the space $B(G)$ of bounded real-valued functions on $G$. The group $G$ acts\footnote{To be more precise, $(\sigma_g)_{g\in G}$ is an anti-action.} on $B(G)$ by $\sigma_t(f)(g):= f(tg), t,g\in G, f\in B(G).$ 
Let $AP(G)$ denote the subspace of  \emph{almost periodic functions}, namely the set of those $f\in B(G)$ for which $\{\sigma_t(f): t\in G\}\subseteq B(G)$ is pre-compact in the sup-norm $\|.\|_\infty$ on $B(G)$.

The following statement is presumably well known to experts.
 %and in fact true in much higher generality,
 However we give  a %relatively self contained
  proof to increase the readability of the paper.
\begin{proposition}\label{ReturnTimesSplitting}
Let $(X,\B,\mu)$ be a Borel probability space, let $(T_g)_{g\in G}$ be a measure preserving $G$-action on $X$, $B\subseteq G, \mu(B)>0$.
Then there exist functions $\phi_c, \phi_{wm}:G\to \R$, where is $\phi_c$ is almost periodic and non-negative such that $\mu(B\cap T_g\me B)=\phi_c(g)+\phi_{wm}(g)$ and
\begin{align}
 m:= \lim_{n\to \infty} \frac1{|F_n|}\sum_{g\in F_n} \phi_{c}(g)&=\lim_{n\to \infty}  \frac{1}{|F_n|} \sum_{g \in F_n} \mu(T_g\me B\cap B)>0,\\ \lim_{n\to \infty} \frac1{|F_n|}\sum_{g\in F_n} |\phi_{wm}(g)|&=0
\end{align}
for any F\o lner sequence $(F_n)_{n \in \N}$.
\end{proposition}

\begin{pf}
Set $\H=L_2(X,\B,\mu)$. Let  $U_g h:= h\circ T_g$, $g\in G$, $h\in \H$ be the induced unitary anti-action of $G$ on $\H$. Pick a F\o lner sequence $(F_n)_{n \in \N}.$ Consider now the following $(U_g)_{g\in G}$-invariant subspaces of $\H$.
\begin{align}
\H_c & = \Big\{f\in \H: \{U_g f:g\in G\} \mbox{ is precompact in the norm topology} \Big\}\\
\H_{wm}&  = \Big\{f\in \H:\frac 1{|F_n|} \sum_{g\in F_n} |\langle U_g f,f'\rangle| \to 0  \mbox { for all } f' \in \H\Big\}.
\end{align}
By \cite[Theorem 1.9]{BeRo88} $\H=\H_c\oplus \H_{wm}$. Since $\H_c$ does not depend on the particular choice of the F\o lner sequence $(F_n)_{n \in \N}$, $\H_{wm}$ doesn't either.
Set $f:= 1_B$ and choose $f_c\in\H_c, f_{wm}\in \H_{wm}$ such that $f= f_c+f_{wm}$. Set
\begin{align}
\phi_c(g):= \langle  U_g f_c, f_c\rangle, \phi_{wm}(g):=\langle U_g f_{wm}, f_{wm}\rangle,\\
  \mu(T_g\me B\cap B)=\langle U_g f,f \rangle=\phi_c(g)+\phi_{wm}(g).
\end{align}
It follows directly from the definition of $\H_{wm}$ that $\lim_{n\to \infty} \frac1{|F_n|}\sum_{g\in F_n} |\phi_{wm}(g)|=0$.  Note that for $t_1, t_2\in G$
\begin{align}
\| \sigma_{t_1} (\phi_c)-\sigma_{t_2} (\phi_c)\|_\infty=
\sup_{g\in G} | \phi_c(t_1 g)-\phi_c(t_2 g)|=\\ \sup_{g\in G} |\langle U_{t_1g} f_c, f_c\rangle-\langle U_{t_2g} f_c, f_c\rangle|=\\
\sup_{g\in G} |\langle U_g(( U_{t_1} - U_{t_2}) (f_c)), f_c\rangle| \leq \|U_{t_1} f_c-U_{t_2} f_c\|_2,
\end{align} hence pre-compactness of $\{U_t f_c:t\in G\}$ implies pre-compactness of $\{\sigma_t (\phi_c):t\in G\}$, thus $ \phi_c$ is almost periodic. By the mean ergodic theorem
\begin{align}
\label{PhiPositive}
 \lim_{n\to \infty} \frac{1}{|F_n|} \sum_{g \in F_n} \phi_c(g)=\lim_{n\to \infty}  \frac{1}{|F_n|} \sum_{g \in F_n} \int_G fU_g f \, d\mu = \int f P f\, d\mu,
\end{align}
where $P$ denotes the projection from $L_2 (\mu)$ onto the subspace of the $U_g$-invariant functions.
Since $\int P f \,d\mu= \int f \, d\mu = \mu(B)$, $f\neq 0$. Thus
\begin{align} 0<\int (Pf)^2\, d\mu= \int Pf P f \, d\mu= \int f P^2 f \, d\mu= \int f P f \, d\mu.\end{align} Hence also the right hand side of (\ref{PhiPositive}) is positive.
\qed\end{pf}

\medskip

We will need the following  alternative characterization of almost periodicity. (See  \cite{BeJH89} for a proof that these two properties are equivalent.)
\begin{lemma}\label{BohrAP}
A function $\phi: G\to \R$ is almost periodic if and only if there exists a continuous function $f:bG\to \R$ such that $\phi=f\circ \iota$.
\end{lemma}

As a consequence of Proposition \ref{ReturnTimesSplitting} and Lemma \ref{BohrAP} we obtain  F\o lner's Theorem (\cite{Foln54a, Foln54b}) for countable amenable groups:
\begin{corollary}\label{AmenableFolner}
Let $G$ be a countable amenable group and let $A\subseteq G$ such that $d^*(A)>0$. Then there exist a Bohr set $B$ and a set $N\subseteq G$ with $d^*(N)=0$ such that
\begin{align}
B\subseteq AA\me\cup N.
\end{align}
\end{corollary}

\begin{pf}
By Lemma \ref{DifferencesContainReturnTimes} there exist a Borel probability space $(X,\B,\mu)$,  $B\in \B,$ $ \mu(B)>0$ and a measure preserving action $(T_g)_{g\in G}$ on $X$, such that $\{g\in G: \mu(T_g\me B\cap B)>0\}\subseteq AA\me$. Pick $m$ and $\phi_c, \phi_{wm} $ according to Proposition \ref{ReturnTimesSplitting} such that $\mu(T_g\me B\cap B)=\phi_c(g)+\phi_{wm}(g)$ for $g\in G$. Set $N=\{g: \phi_{wm}< -m/2\}$ and $\psi= \phi_{c}-m/2$. Then $ d^*(N)=0$ and for $g\in G\setminus N, \psi(g)>0$ implies that $\mu(T_g\me B\cap B)>0$. Pick a continuous function $f:bG\to \R$ such that $\psi=f\circ \iota$.  Since $\lim_{n\to \infty}\frac1{|F_n|}\sum_{g\in F_n} \psi(g)=m/2$,  $f $ takes some positive value, in particular $U:=\{x\in bG: f(x)>0\}$ is a non-empty open set. Putting things together we have
\begin{align}
\iota\me U= \{g:\psi(g)>0\}\subseteq \{g:\mu(T_g\me B\cap B)>0\}\cup N\subseteq AA\me \cup N.
\end{align}
 \qed
 \end{pf}

\medskip

Having Corollary \ref{AmenableFolner} at hand, Theorem \ref{AmenableBohrness} follows from Proposition \ref{PiecewiseDD} once we establish the following regularity property of piecewise Bohr sets.

\begin{lemma}\label{PBRegularity}
Let $S,T\subseteq G$. If  $S$ is piecewise Bohr and for each finite set $ H\subseteq G$ there is some $t_H\in G$ such that $(S\cap H) t_H\subseteq T$ then $T$ is piecewise Bohr as well.
\end{lemma}

\begin{pf}
There exist a thick set $H\subseteq G$ and an open set $ U\subseteq bG$ such that $H\cap \iota\me[U]\subseteq S$. Pick sequences $(H_n)_{n\in \N}$ and $(s_n)_{n\in\N}$ such that $H_n s_n\uparrow G$ and $H_n \subseteq H$. Pick for each $n\in \N$ some $t_n$ such that $(\iota\me [U]\cap H_n)t_n\subseteq T$.
Then
\begin{align}
T\supseteq (\iota\me [U]\cap H_n)t_n =\{g\in H_n:\iota(g)\in U\} t_n=
\{gt_n\in H_nt_n:\iota(g)\in U\}=\\ \{h\in H_n t_n: \iota(h)\iota(t_n\me)\in U\}=\{h\in H_n t_n: \iota(h)\in U\iota(t_n)\}=\iota\me [U\iota(t_n)]\cap H_nt_n
\end{align}
Choose an accumulation point $x $ of $\iota(t_n)\me, n=1,2, \ldots$ and open sets $U_1, U_2$ such that $x\in U_2$ and $U_1\cdot U_2\subseteq U$. Then $U_1 \iota (t_n)\me\subseteq  U$ for infinitely many $n\in \N$ and for each such $n$
\begin{align}
\iota\me [U_1] \cap H_nt_n\subseteq T,
\end{align} hence $T$ is piecewise Bohr.
\qed\end{pf}

\medskip

\newproof{potAB}{Proof of Theorem \ref{AmenableBohrness}}

\begin{potAB}
Pick the set $D$ in $G$ of positive upper Banach density guaranteed by Proposition \ref{PiecewiseDD}. Then by Corollary \ref{AmenableFolner} the set $DD\me$ is piecewise Bohr.  By Lemma \ref{PBRegularity} the set $AB$ is piecewise Bohr.
\qed\end{potAB}

\medskip

\section{Abelian versus non-abelian}
\label{AbelianSection}

The following Lemma is the only remaining fact needed for the proof
of Theorem \ref{cor_pwBohr_char}.

\begin{lemma}\label{ThickIsDivisible}
Assume that $(G,+)$ is a countable  abelian group and $T\subseteq G$ is thick. Then there exist thick sets $T_1,T_2\subseteq G$ such that $T_1+T_2\subseteq T$.
\end{lemma}

\begin{pf}
Pick sequences $(c_n)_{n \in \N}, (K_n)_{n \in \N}$ such that all $K_n\subseteq G$ are finite, $K_n\uparrow G$ and $\bigcup_{n\in\N} K_n+c_n\subseteq T$. We will inductively define sequences $(a_n)_{n \in \N},(b_n)_{n \in \N}$ such that
\begin{align}
\bigcup_{l\in\N} K_l+a_l+\bigcup_{m\in\N} K_m+b_m\subseteq T.
\end{align}
To start the induction, let $a_1\in G$ be arbitrary, pick $n$ such that $K_1+a_1+K_1\subseteq K_n$ and set $b_1=c_n$ such that $K_1+a_1+K_1+b_1\subseteq K_n+c_n\subseteq T$.

Next assume that after $k$ steps $a_1,\ldots , a_k,$ $ b_1,\ldots b_k\in G$ have been chosen such that $ \bigcup_{l\leq k} K_l+a_l+\bigcup_{m\leq k} K_m+b_m\subseteq T$. Pick $n$ such that  $ K_{k+1} + \bigcup_{l\leq k} K_l+a_l+\bigcup_{m\leq k} K_m+b_m\subseteq K_n$ and set $a_{k+1}:=c_n$. Choose $b_{k+1}$ analogously. The induction continues.
\qed\end{pf}

\medskip

\newproof{potLA}{Proof of Theorem \ref{cor_pwBohr_char}}
\begin{potLA}
If $C \subseteq G $ is piecewise Bohr then $C \supseteq B \cap T$,
where $B$ is a Bohr set and $T$ is a thick set. As explained in
Subsection \ref{subsection_piecewise_Bohr} one can find Bohr sets
$B_0,B_1 \subset G $ such that $B \supseteq B_0+B_1$. By Lemma
\ref{ThickIsDivisible} we can find thick sets $T_1,T_2$ in $G$ such
that $T_1+T_2 \subseteq T$. Then $C \supseteq (B_0\cap T_1) +
(B_1\cap T_2)$. On the other hand, if $A+B \subseteq C$ for $A,B$ of
positive upper Banach density then, by Theorem
\ref{AmenableBohrness}, $C$ is piecewise Bohr. \qed\end{potLA}

\medskip
% Property of the Bohr set used in the Proof of Theorem \ref{cor_pwBohr_char} implies the following statement.

One may wonder  whether given three sets $A,B,C$ of positive upper Banach density in an abelian group the sum $A+B+C$ has stronger properties then sumset of two sets. The following result, which follows from the familiar by now fact that a piecewise  Bohr set contains the sum of two piecewise Bohr sets, shows that there is not much to look for.
\begin{proposition}
\label{k_folded_cor} Let $G$ be a countable abelian group and let
$A,B \subseteq G$ be such that $d^{*}(A),d^{*}(B)>0$. Then for every
$k \in \N$ there exist piecewise Bohr sets $C_1,\ldots,C_k$ such
that
\[
C_1+C_2+\ldots+C_k \subseteq A+B.
\]
\end{proposition}

The following Proposition \ref{ThickIsNotDivisible} demonstrates
that in Theorem \ref{cor_pwBohr_char} one cannot drop the assumption
of commutativity of the group $G$. However, before formulating
Proposition \ref{ThickIsNotDivisible} we want to introduce some
convenient terminology. Note first that the definition of upper
Banach density introduced in Subsection %\marginpar{\tiny{In my opinon it is quite misleading to discuss right Banach density here. There is no reason to expect that a thick set should contain the product of two sets of positive right upper Banach density, in fact it is straight forward to find a thick set in the Heisenberg group which doesn't contain a set of right positive upper Banach density. Also it might ``mislead'' the referee to ask more questions on wether $AB$, where $A$ is leftlarge and $B$ is right large contains something. Could we kill the paragraph? Is this font readable to you?}}
\ref{subsection_amenable_groups} is based on the notion of left F\o
lner sequence. One could also introduce a ``right" version of upper
Banach density with the help of the notion of right F\o lner
sequences (that is a sequence satisfying
%\begin{align}
$\lim_{n\to \infty} \frac{|F_n g \bigtriangleup F_n|}{|F_n|}=0
\textnormal{)}$.
%\end{align}
Accordingly, we will say that a set $ A \subseteq G $ is \emph{left
large} (\emph{right large}) if it has positive upper ``left"
(``right") Banach density. Finally, let us say that a set $A
\subseteq G$ is \emph{large} if it is either left large or right
large.

\begin{proposition}\label{ThickIsNotDivisible}
Let $G$ be the Heisenberg group over the integers, i.e.\ the group
of $3\times 3$ upper triangular matrices with integer entries and
1's on the diagonal. There exists a thick set $T\subseteq G$ which
does not contain the product $AB$ of any two large sets
$A,B\subseteq G$.
\end{proposition}

\begin{pf}
We will view $G$ as $\Z^3$ equipped with the operation given by
\begin{align}
\left(a\xx,a\yy,a\zz\right)*\left(b\xx,b\yy,b\zz\right):= \left(a\xx+b\xx, a\yy+b\yy, a\zz+b\zz+a\xx b\yy\right).
\end{align} Set $K_n=\{-n,\ldots, n\}^3$ for $n\in \N$ and $T=\Cup_{n\in \N} K_n * (n^2,0,0)$.  Assume that, contrary to the claim of our
 Proposition, there exist large sets $A,B\subseteq G$ such that  $A*B\subseteq T$. Pick $b_1=(b_1\xx,b_1\yy,b_1\zz),$ $b_2=(b_2\xx,b_2\yy,b_2\zz)\in B$ such that $b_1\yy\neq b_2\yy$. Set $n_0=10 \left(|b_1\xx|+|b_1\yy|+|b_1\zz|+|b_2\xx|+|b_2\yy|+|b_2\zz|\right)$. Since $A$ is infinite, $Ab_1$ is not contained in $\Cup_{n\leq n_0} K_n*(n^2,0,0)$. Hence there exist $a=(a\xx,a\yy,a\zz)\in A$ and $m\geq n_0$ such that $a*b_1\in K_m  * (m^2,0,0)$. Note that this implies that $a\xx\in [m^2-2m,m^2+2m]$.  By assumption, $a*b_2\in T$ and since the difference $\big|\big(a\xx+b_1\xx\big)-\big(a\xx+b_2\xx\big)\big|$ is small compared to $m$, we have in fact $a*b_2\in K_m*(m^2,0,0)$. This implies that the $z$-coordinates of $a*b_1$ and $a*b_2$ differ at most by $2m$, hence
\begin{align}
2m&\geq \left|\left( a\zz+b_1\zz+a\xx b_1\yy\right)-\left( a\zz+b_2\zz+a\xx b_2\yy\right)\right|\\
&=\left| b_1\zz-b_2\zz + a\xx \left( b_1\yy-b_2\yy \right)\right|
\end{align}
which is not possible since $\big|b_1\yy-b_2\yy\big|\geq 1$ and $a\xx$ is of order $m^2$.
\qed\end{pf}

\medskip

\begin{acknowledgment}
The authors thank Michael Hochman, Gabriel Maresch and Ilya Shkredov for helpful comments on the topic of this paper.
\end{acknowledgment}

%\cite{Jin02, JiKe03, Jin04}
%\cite{BeFW06}
%\cite{BeHK05}
%\bibliographystyle{alpha}%{ams-pln}
%\bibliography{MathiasB}

\begin{thebibliography}{BBHS06}

\bibitem[BBHS06]{BBHS06}
M.~Beiglb{\"o}ck, V.~Bergelson, N.~Hindman, and D.~Strauss.
\newblock Multiplicative structures in additively large sets.
\newblock {\em J. Combin. Theory Ser. A}, 113(7):1219--1242, 2006.

\bibitem[Ber85]{Berg85}
V.~Bergelson.
\newblock Sets of recurrence of {${\bf Z}\sp m$}-actions and properties of sets
  of differences in {${\bf Z}\sp m$}.
\newblock {\em J. London Math. Soc. (2)}, 31(2):295--304, 1985.

\bibitem[Ber87]{Ber85}
V.~Bergelson.
\newblock Ergodic {R}amsey theory.
\newblock {\em Amer. J. Math.}, 65(6):63--87, 1987.

\bibitem[Ber00]{Be00}
V.~Bergelson.
\newblock Ergodic theory and {D}iophantine problems.
\newblock In {\em Topics in symbolic dynamics and applications (Temuco, 1997)},
  volume 279 of {\em London Math. Soc. Lecture Note Ser.}, pages 167--205.
  Cambridge Univ. Press, Cambridge, 2000.

\bibitem[Ber06]{Berg06}
V.~Bergelson.
\newblock Combinatorial and {D}iophantine applications of ergodic theory.
\newblock In {\em Handbook of dynamical systems. {V}ol. 1{B}}, pages 745--869.
  Elsevier B. V., Amsterdam, 2006.
\newblock Appendix A by A. Leibman and Appendix B by Anthony Quas and
  M{\'a}t{\'e} Wierdl.

\bibitem[BFW06]{BeFW06}
V.~Bergelson, H.~Furstenberg, and B.~Weiss.
\newblock Piecewise-{B}ohr sets of integers and combinatorial number theory.
\newblock In {\em Topics in discrete mathematics}, volume~26 of {\em Algorithms
  Combin.}, pages 13--37. Springer, Berlin, 2006.

\bibitem[BJM89]{BeJH89}
J.~Berglund, H.~Junghenn, and P.~Milnes.
\newblock {\em Analysis on semigroups}.
\newblock Canadian Mathematical Society Series of Monographs and Advanced
  Texts. John Wiley \& Sons Inc., New York, 1989.
\newblock Function spaces, compactifications, representations, A
  Wiley-Interscience Publication.

\bibitem[BM98]{BeMc98}
V.~Bergelson and R.~McCutcheon.
\newblock Recurrence for semigroup actions and a non-commutative {S}chur
  theorem.
\newblock In {\em Topological dynamics and applications (Minneapolis, MN,
  1995)}, volume 215 of {\em Contemp. Math.}, pages 205--222. Amer. Math. Soc.,
  Providence, RI, 1998.

\bibitem[BR88]{BeRo88}
V.~Bergelson and J.~Rosenblatt.
\newblock Mixing actions of groups.
\newblock {\em Illinois J. Math.}, 32(1):65--80, 1988.

\bibitem[F{\o}l54a]{Foln54a}
E.~F{\o}lner.
\newblock Generalization of a theorem of {B}ogolio\`uboff to topological
  abelian groups. {W}ith an appendix on {B}anach mean values in non-abelian
  groups.
\newblock {\em Math. Scand.}, 2:5--18, 1954.

\bibitem[F{\o}l54b]{Foln54b}
E.~F{\o}lner.
\newblock Note on a generalization of a theorem of {B}ogolio\`uboff.
\newblock {\em Math. Scand.}, 2:224--226, 1954.

\bibitem[Fur81]{Furs81}
H.~Furstenberg.
\newblock {\em Recurrence in ergodic theory and combinatorial number theory}.
\newblock Princeton University Press, Princeton, N.J., 1981.
\newblock M. B. Porter Lectures.

\bibitem[Hin74]{Hind74}
N.~Hindman.
\newblock Finite sums from sequences within cells of a partition of {$N$}.
\newblock {\em J. Combinatorial Theory Ser. A}, 17:1--11, 1974.

\bibitem[Hin82]{Hind82}
N.~Hindman.
\newblock On density, translates, and pairwise sums of integers.
\newblock {\em J. Combin. Theory Ser. A}, 33(2):147--157, 1982.

\bibitem[HS98]{HiSt98}
N.~Hindman and D.~Strauss.
\newblock {\em Algebra in the {S}tone-\v {C}ech compactification}, volume~27 of
  {\em de Gruyter Expositions in Mathematics}.
\newblock Walter de Gruyter \& Co., Berlin, 1998.
\newblock Theory and applications.

\bibitem[Jin02]{Jin02}
R.~Jin.
\newblock The sumset phenomenon.
\newblock {\em Proc. Amer. Math. Soc.}, 130(3):855--861, 2002.

\bibitem[Jin04]{Jin04}
R.~Jin.
\newblock Standardizing nonstandard methods for upper {B}anach density
  problems.
\newblock In {\em Unusual applications of number theory}, volume~64 of {\em
  DIMACS Ser. Discrete Math. Theoret. Comput. Sci.}, pages 109--124. Amer.
  Math. Soc., Providence, RI, 2004.

\bibitem[Jin08]{Ji08}
R.~Jin.
\newblock {\em Private Communication}, 2008.

\bibitem[JK03]{JiKe03}
R.~Jin and H.~Keisler.
\newblock Abelian groups with layered tiles and the sumset phenomenon.
\newblock {\em Trans. Amer. Math. Soc.}, 355(1):79--97, 2003.

\bibitem[Lin01]{Lind01}
E.~Lindenstrauss.
\newblock Pointwise theorems for amenable groups.
\newblock {\em Invent. Math.}, 146(2):259--295, 2001.

\bibitem[OW87]{OrWe87}
D.~Ornstein and B.~Weiss.
\newblock Entropy and isomorphism theorems for actions of amenable groups.
\newblock {\em J. Analyse Math.}, 48:1--141, 1987.

\bibitem[Wei01]{Weis01}
B.~Weiss.
\newblock Monotileable amenable groups.
\newblock In {\em Topology, ergodic theory, real algebraic geometry}, volume
  202 of {\em Amer. Math. Soc. Transl. Ser. 2}, pages 257--262. Amer. Math.
  Soc., Providence, RI, 2001.

\end{thebibliography}

\def\ocirc#1{\ifmmode\setbox0=\hbox{$#1$}\dimen0=\ht0 \advance\dimen0
  by1pt\rlap{\hbox to\wd0{\hss\raise\dimen0
  \hbox{\hskip.2em$\scriptscriptstyle\circ$}\hss}}#1\else {\accent"17 #1}\fi}

\end{document}